\documentclass[12pt,a4paper]{article}

\usepackage{amsmath}
\usepackage{amssymb}
\usepackage{amsthm}
\usepackage{amsfonts}
\usepackage{graphicx}
\usepackage[all]{xy}
\usepackage{verbatim}
\usepackage{parskip}
\usepackage{csquotes}
\usepackage{tkz-graph}
\usepackage[left=2cm,right=2cm,top=3cm,bottom=2.5cm]{geometry}
\usepackage{soul}

\usepackage{hyperref}
\usepackage{ulem}

\newcommand{\norm}[2]{\left\| #1 \right\|_{#2}}
\newcommand{\M}{\mathcal{M}}

\newcommand{\R}{\mathbb{R}}

\theoremstyle{definition}
\newtheorem{theorem}{Theorem}[section]
\newtheorem{definition}[theorem]{Definition}

\newtheorem{remark}[theorem]{Remark}
\newtheorem{proposition}[theorem]{Proposition}

\newtheorem{lemma}[theorem]{Lemma}

\begin{document}
\thispagestyle{empty}
\begin{center}
    {\LARGE Dynamic Interpolation and Higher-Order Turbulence Scales in 3D Navier--Stokes Flows \par}
    
    \vspace{0.8em}
    
{\large
    Moses Patson Phiri \qquad Zoran Gruji\'c
    \par
}
 
    \vspace{0.2cm}

    {\normalsize
    Department of Mathematics,
    University of Alabama at Birmingham
    \par}

    \vspace{0.4cm}

    {\normalsize \today\par}
\end{center}

\vspace{0.4cm}
\vspace{0.8cm}

\begin{abstract}
We investigate the three-dimensional unforced Navier–Stokes equations through the interaction between spatial analyticity and the geometric sparseness of regions of intense fluid activity (i.e., regions where the fluid’s motion varies sharply). Interpreting higher-order derivative ratios as turbulence scales, we introduce time-weighted bridge inequalities that quantify the relationship between lower- and higher-order derivative norms. The time weights give the derivative ratios greater flexibility while ensuring that the resulting lower bound on the radius of analyticity dominates the sparseness scale of the super-level sets, indicating that the flow has entered the dissipation regime. Combining this scale comparison with a harmonic-measure maximum principle on one-dimensional sparse sets yields a contraction in the relevant derivative norm, thereby excluding a class of strongly turbulent blow-up scenarios. Consequently, the strong solution extends analytically beyond the candidate blow-up time.
\end{abstract}

\textbf{Keywords: }3D Navier–Stokes equations; spatial analyticity; analyticity radius; geometric sparseness; dynamic interpolation; finite-time blow-up 

 \newpage 

\section{Introduction}\label{section1}

Whether the internal competition between nonlinear transport and viscous dissipation in the three-dimensional unforced Navier--Stokes equations (NSE) can generate finite-time singularities from smooth initial data remains a fundamental open problem in the mathematical theory of fluid flow \cite{Fefferman2000,RobinsonRodrigoSadowski2016,Tsai2018}. The NSE are the fundamental equations governing the
motion of 3D incompressible, viscous, Newtonian fluids and are given by
\begin{align}
\label{NSE}
\partial_t u + (u\cdot\nabla)u &= \nu\Delta u-\nabla p,
&\text{ in } \R^3 \times (0,T),\\
\nabla\cdot u &= 0,
&\text{ in } \R^3\times(0,T),\\
u(\cdot,0) &= u_0(\cdot),
&\text{ in } \R^3\times\{t=0\}.
\end{align}
The vector field $u$ and the scalar field $p$ denote the fluid velocity and pressure, respectively. Throughout our analysis, the viscosity is normalized
to $\nu=1$, and we consider the unforced case $f=0$.

The NSE possesses a natural scaling that leaves the system invariant. Let
$\lambda>0$. If $(u,p)$ is a solution of the NSE, then
\[
u_\lambda(x,t)
=
\lambda u(\lambda x,\lambda^2 t),
\qquad
p_\lambda(x,t)
=
\lambda^2 p(\lambda x,\lambda^2 t)
\]
is also a solution. This scaling structure is closely tied to questions of regularity and possible finite-time singularity formation. For the three-dimensional unforced system considered here, it remains unknown whether smooth initial data can develop a singularity in finite time. A central difficulty arises from the competition between nonlinear transport, which can
transfer activity toward progressively smaller spatial scales, and viscous diffusion, which acts to dissipate such concentrations
\cite{ConstantinFoias1988,FoiasManleyRosaTemam2001}.

A major difficulty is the scaling gap between the available a priori bounds
and known regularity criteria. The Leray–Hopf energy bounds yield
\[u\in L^{p}(0,T;L^{q}(\mathbb{R}^{3})),
\qquad
\frac{2}{p}+\frac{3}{q}=\frac{3}{2},
\qquad 2\le q\le 6,
\]
whereas the \text{Ladyzhenskaya--Prodi--Serrin} (LPS) criterion guarantees regularity under the stronger integrability condition

\[
u\in L^{p}(0,T;L^{q}(\mathbb{R}^{3})),
\qquad
\frac{2}{p}+\frac{3}{q}=1,
\qquad 3<q\le\infty.
\]

Thus, the energy estimates alone do not reach the scaling-critical LPS regime, leaving the gap
\[
\dfrac{3}{2}\longrightarrow 1.
\]
This motivates looking beyond purely integrability-based criteria and
exploiting geometric intermittency, as developed by Gruji\'c and co-authors
\cite{grujicz,grujic2,Bradshaw,GrujicKukavica1998}, in particular the
sparseness of regions of intense fluid activity. Philosophically, the emphasis shifts from controlling the size of the solution globally to controlling the geometry of the regions where the solution is large.

This approach combines spatial analyticity and its associated higher-order derivative estimates \cite{Guberovic2010} with the spatial sparseness of derivative super-level sets. When the analyticity radius dominates the sparseness scale, harmonic-measure estimates and the maximum principle yield the norm control needed to exclude blow-up \cite{GrujicXu2024,grujic}. Albritton and Bradshaw developed alternative sparseness-based regularity criteria and highlighted limitations in interpreting higher-order sparseness estimates alone as a reduction of the scaling gap \cite{AlbrittonBradshaw}. Here, additional time-weighted relations between derivative levels ensure this scale comparison under the hypotheses of our main theorem.

    In this work, we investigate this interaction at the level of higher-order spatial derivatives and derive a no-blow-up extension criterion for strong solutions of the three-dimensional Navier--Stokes equations. The mechanism is based on a time-weighted comparison across derivative levels, which enables a
\textit{dynamic interpolation} between higher-order scales. For $2\ell\le q\le k$, we introduce the modified time-weighted bridge
\[
\|D^q u(t)\|_\infty^{\tfrac{1}{q+1}}
\le
\mathcal{M}_{q,k}\,
\|D^k u(t)\|_\infty^{\tfrac{1}{k+1}},
\]
where
\begin{align}
\label{twedge}
\mathcal{M}_{q,k}
&=
\tau_k^{\beta\left(
\tfrac{1}{(q+1)^\alpha}
-
\tfrac{1}{(k+1)^\alpha}
\right)}
\tfrac{(q!)^{\tfrac{1}{q+1}}}
     {(k!)^{\tfrac{1}{k+1}}},
\qquad
0<\alpha\le1,\quad \beta>0,
\end{align}
and $\tau_k\in(0,1)$ denotes the associated time weight. The factor
$\mathcal{M}_{q,k}$ transfers information from the highest derivative level $k$ to the intermediate levels $q<k$ while retaining the time dependence of
the evolution. Here, $\alpha$ determines how the temporal weight is
distributed across the derivative hierarchy, whereas $\beta$ controls the strength of the resulting gain. The bridge therefore provides a dynamic interpolation across higher-order derivative levels and the flexibility needed
to compare the analyticity and sparseness scales as a potential singular time is approached.

\begin{remark}[Sign of the time-weight exponent] In Phiri's hyperdissipative Navier--Stokes analysis \cite{Phiri}, the time-weight exponent is negative: \[ -\frac{1}{2\beta_{\mathrm{HD}}} \left( \frac{\mu+j}{j+1} - \frac{\mu+k}{k+1} \right) = -\frac{\mu-1}{2\beta_{\mathrm{HD}}} \left( \frac{1}{j+1} - \frac{1}{k+1} \right) <0, \qquad j<k, \qquad \mu>1. \] Consequently, the time weight grows as the remaining time approaches zero, permitting growth in the derivative ratios allowed by the bridge. In the present construction for the classical Navier--Stokes equations, the modified bridge has a positive time-weight exponent: \[ \beta \left( \frac{1}{(q+1)^{\alpha}} - \frac{1}{(k+1)^{\alpha}} \right) >0, \qquad q<k, \qquad 0<\alpha\le 1, \qquad \beta>0. \] Because $0<\tau_k<1$, this weight strengthens the derivative comparison and supplies an additional smallness factor. Under the hypotheses of the main theorem, this gain is sufficient to close the analyticity--sparseness gap. Thus, while the hyperdissipative construction uses time weights to relax the derivative-ratio bounds, the present construction uses them to sharpen those bounds. Here, $\beta_{\mathrm{HD}}>1$ denotes the hyperdissipation exponent and is distinct from the time-weight parameter $\beta$. \end{remark}

Large values of $\|D^k u(t)\|_\infty$ correspond to increasingly small
spatial scales, motivating the interpretation of higher-order derivatives as a hierarchy of turbulence scales. For $0<\lambda<1$, define the associated superlevel set by
\[
S_{\lambda}^{(k)}(t)
:=
\left\{
x\in\R^3:
|D^k u(x,t)|
>
\lambda\|D^k u(t)\|_\infty
\right\}.
\]
The spatial sparseness of $S_{\lambda}^{(k)}(t)$ provides the geometric scale to be compared with the analyticity radius.

\begin{theorem}[Exclusion of finite-time blow-up]
\label{Mainthm2}
Let $\alpha\in(0,1)$ and $\beta>0$ be arbitrary constants. Let $T^{*} < \infty$ be the
first possible blow-up time, and assume the hypotheses of
Theorem~\ref{MainThm1}. Assume that there exists \(\varepsilon > 0\) and an integer \(K_{0} \geq 2\ell\) such that, for every \(k\geq K_{0},\) the function
\[t \longmapsto \norm{D^ku(t)}{\infty} \]
is strictly increasing on the interval \((T^{*} - \varepsilon, T^{*})\). If
\begin{equation}
\label{requiredthresh}
\alpha\beta(\ell+1)^{1-\alpha}>\dfrac{1}{4},
\end{equation}
then $T^*$ is not a blow-up time. Consequently, the solution $u(\cdot,t)$
extends spatially analytically beyond $T^*$.
\end{theorem}

\begin{remark}
The index $\ell$ in Theorem~\ref{Mainthm2} is the same index furnished by Theorem~\ref{MainThm1}; hence, no additional largeness assumption on $\ell$ is imposed. The threshold \eqref{requiredthresh} is precisely the condition required for the analyticity radius to dominate the corresponding sparseness scale.
\end{remark}

\begin{remark}
The parameters $\alpha$ and $\beta$ may be chosen in two natural ways.

\begin{enumerate}
    \item[(i)] For any fixed $\alpha\in(0,1)$, choosing
    \[
    \beta>\frac{1}{4\alpha}
    \]
    guarantees that \eqref{requiredthresh} holds.

    \item[(ii)] More generally, for any $\alpha\in(0,1)$ and $\beta>0$, one may choose $\ell$ sufficiently large so that
    \eqref{requiredthresh} holds, provided that this choice is compatible with the admissibility condition in \eqref{rangeoindices},
    \[
    (\ell!)^{2} \leq \norm{u_0}{\infty}.
    \]
\end{enumerate}

These choices are sufficient to close the contradiction argument in the proof of Theorem~\ref{Mainthm2}. In particular, $\alpha=0.89$ is admissible for suitable $\beta$ and $\ell$. This value is consistent with the computational results of Gruji\'c and Mohebujjaman \cite{GrujicMohebujjaman2026}, where
higher-order derivative ratios for the three-dimensional Taylor--Green vortex were observed to obey a temporal scaling law with $\alpha\approx0.89$ near the
enstrophy peak.
\end{remark}

The mechanism behind Theorem~\ref{Mainthm2} becomes apparent in the Duhamel
formulation of the Navier--Stokes equations, with the underlying computation given explicitly in \cite{GrujicXu2024}. After applying $D^k$ to the integral equation, the nonlinear term generates products involving several intermediate
derivative levels. In the classical approach, these lower-order factors are controlled through Gagliardo--Nirenberg interpolation. Here, instead, we apply the derivative bridge at the relevant levels, transferring control from the
$k$th derivative to the intermediate derivatives while retaining quantitative information on the time scale of the evolution. The resulting interpolation is therefore dynamic, coupling the derivative hierarchy to the temporal evolution of the flow.

The effect of this interpolation is reflected directly in the corresponding lower bounds on the analyticity radius. At the level of scaling, the classical
Gagliardo--Nirenberg interpolation yields
\[
\rho_k
\gtrsim
\|D^k u(t)\|_\infty^{
-\frac{3}{2}\frac{1}{k+\frac32}}.
\]
The time-independent derivative bridge improves this estimate to
\[
\rho_k
\gtrsim
\|D^k u(t)\|_\infty^{-\frac{1}{k+1}}.
\]
However, the corresponding geometric sparseness scale is of order
\[
r_k
\sim
\|D^k u(t)\|_\infty^{-\frac{1}{k+\frac{3}{2}}},
\]
and
\[
\frac{1}{k+1}
>
\frac{1}{k+\tfrac{3}{2}}.
\]

Although the time-independent bridge improves the analyticity radius, a scaling gap remains. Our construction introduces time weights that provide additional flexibility in the derivative ratios, yielding a stronger lower bound on the analyticity radius sufficient to dominate the geometric sparseness scale. This allows us to exclude a broader class of blow-up scenarios.

The time-weighted derivative bridge provides the additional gain required to close this gap. More precisely, the resulting analyticity estimate takes the form
\[
\rho_k
\gtrsim
\|D^k u(t)\|_\infty^{
-\tfrac{1}{
k+1+2\alpha\beta(\ell+1)^{1-\alpha}}}.
\]
Consequently, since
\[
\alpha\beta(\ell+1)^{1-\alpha}>\frac14
\]
we have that
\[
k+1+2\alpha\beta(\ell+1)^{1-\alpha}
>
k+\frac32,
\]
and hence
\[
\frac{1}{
k+1+2\alpha\beta(\ell+1)^{1-\alpha}}
<
\frac{1}{k+\tfrac{3}{2}}.
\]
In the large-derivative regime, the resulting lower bound on the analyticity radius therefore reaches and eventually dominates the sparseness scale. This
scale separation allows the harmonic-measure maximum principle to yield a
contraction at a suitable escape time, closing the contradiction argument and establishing the no-blow-up extension criterion.

The principal new ingredient is therefore the time-weighted derivative
bridge. While the time-independent bridge already improves the classical
interpolation across derivative levels, the additional time weight couples this comparison to the temporal evolution and supplies the gain needed to cross the sparseness threshold. The higher-order derivative ratios may consequently
be viewed as a hierarchy of dynamically coupled spatial scales, linking the interpolation mechanism to their interpretation as higher-order turbulence scales. In this way, the analyticity--sparseness framework is extended by a time-dependent interpolation mechanism that converts higher-order derivative information quantitatively into a regularity mechanism.

The remainder of the paper is organized as follows. Section~\ref{prelim}
collects the analyticity, sparseness, and harmonic-measure preliminaries used throughout the paper. Section~\ref{results} develops the time-weighted derivative bridge and the associated higher-order estimates, compares the
resulting analyticity scales, establishes the analyticity--sparseness gap, and combines this scale separation with the harmonic-measure argument to prove the main regularity result.

\newpage
\section{Preliminaries}\label{prelim}

We begin by recalling the notions and baseline facts concerning the
\textit{sparseness of regions of intense fluid activity}, spatial analyticity,
and harmonic measure that will be used throughout the paper.

\subsection{Higher-order derivatives and maximal existence time}

We use the standard local existence and continuation theory
for strong solutions of the Navier--Stokes equations; see
\cite{RobinsonRodrigoSadowski2016,Tsai2018}. Throughout, $u$ denotes a strong solution of the 3D
unforced NSE on $[0,T^*)$, where $T^*>0$ denotes the
maximal time of existence. 
For $k\ge0$, $D^k u : = \nabla^k u$ denotes the collection of all spatial derivatives
$D^\gamma u_j$ with $|\gamma|=k$ and $1\le j\le3$. Throughout the paper,
we use the convention
\[
\|D^k u(t)\|_\infty
:=
\max_{\substack{|\gamma|=k\\1\le j\le3}}
\|D^\gamma u_j(\cdot,t)\|_{L^\infty(\mathbb R^3)}.
\]
If $T^*<\infty$, we refer to $T^*$ as the finite-time
blow-up time.

\subsection{Sparseness of regions of intense activity}

We recall the notion of directional sparseness.

\begin{definition}[1D $\delta$-sparseness]
\label{def:1D-sparse}
Let $S\subset\mathbb R^3$ be open, $x_0\in\mathbb R^3$,
$\delta\in(0,1)$, and $r>0$. We say that $S$ is
\textit{1D $\delta$-sparse around $x_0$ at scale $r$} if there exists
a unit vector $b$ such that
\[
\frac{
\left|
S\cap(x_0-rb,x_0+rb)
\right|
}{2r}
\le\delta,
\]
where $|\cdot|$ denotes one-dimensional Lebesgue measure along the segment.
\end{definition}

We also recall the corresponding volumetric notion.

\begin{definition}[3D $\delta$-sparseness]
\label{def:3D-sparse}
Let $S\subset\mathbb R^3$ be open. We say that $S$ is
\textit{3D $\delta$-sparse around $x_0$ at scale $r$} if
\[
\frac{|S\cap B_r(x_0)|}{|B_r(x_0)|}
\le\delta.
\]
\end{definition}

In three dimensions, volumetric $\delta$-sparseness at scale $r$ implies
1D $\delta^{1/3}$-sparseness at the same scale in at least one direction
through $x_0$. The converse does not hold. Hence the 1D condition is the
weaker geometric requirement and is the form naturally suited to the
harmonic-measure argument. 

At derivative level $k$, the a priori sparseness scale associated with the
superlevel sets above is of the order
\[
r_k(t)
\approx
\|D^k u(t)\|_\infty^{-\frac{2}{2k+3}}
=
\|D^k u(t)\|_\infty^{-\frac{1}{k+\frac32}}.
\]
We refer to $r_k(t)$ as the \textit{sparseness scale} at derivative level $k$.
The precise multiplicative constants are inherited from the corresponding
geometric sparseness estimate.

\subsection{Spatial analyticity}

Spatial analyticity provides the analytic length scale used in the regularity
argument. For $\rho>0$, define
\[
\mathcal T_\rho
:=
\left\{
x+iy\in\mathbb C^3:
x,y\in\mathbb R^3,\ |y|<\rho
\right\}.
\]

\begin{definition}[Analyticity radius]
For $k\ge0$, we denote by $\rho_k(t)$ the radius of spatial analyticity of
$D^k u(\cdot,t)$, namely
\[
\rho_k(t)
:=
\sup\left\{
\rho>0:
D^k u(\cdot,t)
\text{ admits a bounded holomorphic extension to }
\mathcal T_\rho
\right\}.
\]
\end{definition}

We also recall the Gagliardo--Nirenberg interpolation inequality, which will
be used in the higher-order analyticity estimates.

\begin{lemma}[Gagliardo--Nirenberg \cite{Gagliardo,Nirenberg}]
\label{lem:GN}
Suppose that $p,q,r\in[1,\infty]$, $m,j,d\in\mathbb N$ with
$0\le j<m$, and $\theta\in[j/m,1]$ satisfy
\[
\frac1p
=
\frac{j}{d}
+
\theta\left(
\frac1r-\frac{m}{d}
\right)
+
\frac{1-\theta}{q}.
\]
Then there exists a constant
\[
C=C(m,d,j,q,r,\theta)>0
\]
such that
\[
\|D^j f\|_{L^p}
\le
C\,
\|D^m f\|_{L^r}^{\theta}
\|f\|_{L^q}^{1-\theta}.
\]
\end{lemma}

The following theorem of Gruji\'c and Xu \cite{grujic}, stated for the Navier-Stokes system, shows that ascending derivative chains yield improved analyticity estimates.

\begin{theorem}[Grujic and Xu \cite{grujic}, Fahart and Grujic \cite{farhat}] \label{theorem1}
Let $\beta>1$, $u_0\in L^2$, and $D^{i}u_0\in L^\infty$ for $0\le i\le k$, and suppose that
\begin{equation}\label{notime}
  \norm{D^{q}u_0}{\infty}^{\tfrac{1}{q+1}}
  \;\leq\;
  \M_{q,k}\,\norm{D^{k}u_0}{\infty}^{\tfrac{1}{k+1}}
  \qquad \forall\, \ell \le q \le k .
\end{equation}
Assume further that the constants $\{\M_{j,k}\}$ and the indices $\ell$ and $k$ satisfy
\begin{equation}\label{cond1}
  \sum_{\ell \le i \le j-\ell}
  \binom{j}{i}\,
  \M_{i,k}^{\,i+1}\,
  \M_{j-i,k}^{\,j-i+1}
  \;\lesssim\;  1,
  \qquad \forall\,\, 2\ell \le j \le k ,
\end{equation}
and
\begin{equation}\label{cond2}
  \norm{u_0}{2} \!
  \sum_{0 \le i \le \ell}
  \binom{j}{i}
  \M_{\ell,k}^{\,\tfrac{(\ell+1)(i+\tfrac{3}{2})}{\ell+\tfrac{3}{2}}}
  \M_{j-i,k}^{\,j-i+1}
  \Big((k!)^{\frac{1}{k+1}}\norm{u_0}{\infty}\Big)^{\tfrac{\left(\tfrac{3}{2}-1\right)(\ell-i)}{\ell+\tfrac{3}{2}}}
  \lesssim 1,
\end{equation}
 for all  \( 2\ell \le j \le k, \) for some functions \(\phi\) and \(\psi\). If
\begin{equation}\label{eq:3.4}
  T_{j}^{\tfrac{1}{2}} \;\lesssim\;
  \norm{D^{k}u_0}{\infty}^{-\tfrac{1}{(k+1)}},
\end{equation}
then for any $\ell \le j \le k$ the complexified solution satisfies the upper bound
\begin{equation}\label{eq:3.5}
  \sup_{t\in(0,T_j)} \;\sup_{y\in \mathcal{D}_t} \norm{D^{j}u(\cdot,y,t)}{\infty}
  \;+\;
  \sup_{t\in(0,T_j)} \;\sup_{y\in \mathcal{D}_t} \norm{D^{j}v(\cdot,y,t)}{\infty}
  \;\lesssim\;
  \norm{D^{j}u_0}{\infty} + \norm{D^{k}u_0}{\infty}^{\frac{j+1}{k+1}},
\end{equation}
where \[\mathcal{D}_{t} = : \{x + iy \in \mathbb{C}^3 : |y| \le ct^{\frac{1}{2\beta}}\}.\]
\end{theorem}

In Theorem~\ref{theorem1}, the coefficient $\mathcal{M}_{m,k}$ is given by the definition in \eqref{twedge}.

Theorem~\ref{scalebound} provides the a prior scale  of  sparseness for the super-level sets of the \(k\)-the order derivatives of the velocity field.

\begin{theorem}[Gruji\'c--Xu \cite{grujic}, sparseness-scale estimate] \label{scalebound}
	Suppose \(u\)  is a Leray solution on \((0,T^*)\), with
	\[
	u\in C((0,T^*);L^\infty(\mathbb{R}^3)).
	\]
	Let \(\gamma\in\mathbb{N}^3\) satisfy \(|\gamma|=k\). For each \(t\in(0,T^*)\), then the super-level sets
	\[
	S_{\gamma,\lambda}^{i,\pm}(t)
	:=
	\left\{
	x\in\mathbb{R}^3 : (D^\gamma u)_i^\pm(x,t)>\lambda \|D^\gamma u(t)\|_{L^\infty}
	\right\},
	\qquad 1\le i\le 3
	\]
	 are \(3D\) \(\delta\)-sparse around any spatial point \(x_0\) at scale
	\[
	r^*(t)=c(\|u_0\|_{L^2})\,\|D^\gamma u(t)\|_{L^\infty}^{-\tfrac{1}{k+ \tfrac{3}{2}}},
	\]
	whenever \(r^*(t)\in(0,1]\).
\end{theorem}

\subsection{Duhamel formulation and derivative interactions}

Let $\mathbb P$ denote the Leray projection onto divergence-free vector
fields. Since $\nu=1$ and $f=0$, the Navier--Stokes equations admit the mild
representation \cite{Kato1984}
\[
u(t)
=
e^{(t-t_0)\Delta}u(t_0)
-
\int_{t_0}^{t}
e^{(t-s)\Delta}
\mathbb P\nabla\cdot(u\otimes u)(s)\,ds,
\qquad
0\le t_0<t<T^*.
\]

Applying $D^\gamma$, with $|\gamma|=k$, gives
\[
D^\gamma u(t)
=
e^{(t-t_0)\Delta}D^\gamma u(t_0)
-
\int_{t_0}^{t}
e^{(t-s)\Delta}
\mathbb P D^\gamma\nabla\cdot(u\otimes u)(s)\,ds.
\]
The Leibniz expansion of the nonlinear term produces interactions between
different derivative levels. Precisely,
\[
D^\gamma\big((u\cdot\nabla)u\big)
=
\sum_{\eta\le\gamma}
\binom{\gamma}{\eta}
\big(D^\eta u\cdot\nabla\big)
D^{\gamma-\eta}u.
\]
Thus, the evolution at derivative level $\gamma$ naturally couples the highest
derivative to a family of intermediate derivative levels. These intermediate derivative factors can be estimated using
the classical Gagliardo--Nirenberg interpolation inequalities
\cite{Nirenberg1959}. In the present work, we instead apply
the derivative bridge at the relevant intermediate levels.

\subsection{Escape times}
The contradiction argument is carried out at times when the $k$th-order derivative norm is about to increase further.

\begin{definition}[Escape time]
\label{def:escape-time}
Fix $k\ge0$. A time $t\in(0,T^*)$ is called an
\textit{escape time at derivative level $k$} if
\[
\|D^k u(s)\|_\infty
>
\|D^k u(t)\|_\infty
\qquad
\text{for every }s\in(t,T^*).
\]
\end{definition}

Escape times therefore provide natural reference times from which the
analyticity estimate is propagated forward and compared with the geometry of
the intense regions at a later time.

\subsection{Harmonic measure}
Let $\Omega\subset\mathbb C$ be a domain and let
$K\subset\partial\Omega$ be a Borel set. For $z\in\Omega$, the quantity
\[
\omega(z,\Omega,K)
\]
denotes the harmonic measure of $K$ relative to $\Omega$, evaluated at $z$.
\begin{proposition}[Ransford \cite{Ransford}]
\label{prop:harmonic-measure}
Let $\Omega\subset\mathbb C$ be an open, connected set whose boundary has
nonzero Hausdorff dimension, and let $K$ be a Borel subset of
$\partial\Omega$. Suppose that $u$ is subharmonic on $\Omega$ and satisfies
\[
u(z)\le M,
\qquad z\in\Omega,
\]
and
\[
\limsup_{z\to\xi}u(z)\le m,
\qquad \xi\in K.
\]
Then
\[
u(z)
\le
m\,h(z,\Omega,K)
+
M\bigl(1-h(z,\Omega,K)\bigr),
\qquad z\in\Omega,
\]
where $h(z,\Omega,K)$ denotes the harmonic measure of $K$ relative to
$\Omega$, evaluated at $z$.
\end{proposition}

The way this principle enters the proof is as follows. Fix $x_0\in\mathbb R^3$ and choose a signed component $D^\gamma u_i^\sigma$ that realizes the relevant derivative amplitude at
$x_0$. If the corresponding superlevel set is 1D sparse around $x_0$, let $\nu\in\mathbb S^2$ denote the direction supplied by
Definition~\ref{def:1D-sparse}. Spatial analyticity then allows us to consider the one-dimensional complexification
\[
F(z)
=
D^\gamma u_i^\sigma(x_0+z\nu,s),
\qquad
|z|<\rho.
\]

After translating $x_0$ to the origin, let $K\subset[-\rho,\rho]$ denote
the good portion of the real diameter, namely the complement of the
corresponding superlevel set along the selected direction. The relevant
domain is the slit disk
\[
\Omega
=
\mathbb D_\rho\setminus K.
\]
Since $x_0$ belongs to the selected superlevel set, its translated point
cannot belong to the good set. Hence
\[
0\notin K,
\qquad
0\in\Omega.
\]
On $K$, the superlevel-set definition provides the improved boundary bound,
while spatial analyticity gives the ambient bound in $\mathbb D_\rho$. The sparseness assumption, together with Solynin's
harmonic-measure estimate \cite{Solynin1999}, provides
a quantitative lower bound on $\omega(0,\Omega,K)$.
Proposition~\ref{prop:harmonic-measure} then yields
the contraction used in the no-blow-up argument.

The preceding ingredients constitute the geometric--analytic framework of
the proof. Section~\ref{results} introduces the time-weighted derivative
bridge, quantifies the resulting improvement of the analyticity radius, and
shows that under the threshold condition of the main theorem the improved
analytic scale reaches the sparseness scale.

\newpage

\section{Main Results}
\label{results}

We now introduce the time-weighted derivative bridge that drives the main
estimates of this work. For $2\ell\le q\le k$, define
\begin{align}
\label{defininitionMjk}
\mathcal{M}_{q,k}
:=
\left(\tau_k\right)^{
\beta\left(
\frac{1}{(q+1)^\alpha}
-
\frac{1}{(k+1)^\alpha}
\right)}
\frac{(q!)^{\frac{1}{q+1}}}
     {(k!)^{\frac{1}{k+1}}},
\end{align}
where
\[
0<\alpha\le1,
\qquad
\beta>0,
\qquad
0<\tau_k<1.
\]
With this notation, the time-weighted bridge takes the form
\begin{align}
\label{timeweightedbridge}
\|D^q u(t)\|_\infty^{\frac{1}{q+1}}
\le
\mathcal{M}_{q,k}\,
\|D^k u(t)\|_\infty^{\frac{1}{k+1}},
\qquad
2\ell\le q\le k.
\end{align}

The coefficient \(\mathcal{M}_{q,k}\) combines the factorial scaling across
derivative levels with a time-dependent gain through $\tau_k$. This additional
factor is the key modification of the derivative bridge and will be used to
improve the corresponding analyticity estimates. We first establish the
estimates needed to propagate this bridge across the derivative hierarchy and
then compare the resulting analyticity scale with the geometric sparseness
scale.

\begin{theorem}\label{MainThm1}
    Suppose \(0 <\alpha < 1,\) \(\beta >0,\) and \(\tau_k \in (0, 1)\). Let \(u_0 \in L^2(\R^3)\), \(D^iu_{0} \in L^{\infty}(\R^3)\), for \(0\leq i \leq k\), where \(k\) is a positive integer. Suppose 
    \begin{align}
        \label{bridge}
    \norm{D^qu_0}{}^{\frac{1}{q+1}} \le (\tau_k)^{ \beta\left(\tfrac{1}{(q+ 1)^{\alpha}} - \tfrac{1}{(k + 1)^{\alpha}}\right)} \dfrac{\left(q!\right)^{\frac{1}{q+1}}}{\left(k!\right)^{\tfrac{1}{k+1}}}  \norm{D^ku_0}{}^{\tfrac{1}{k+1}}, \qquad 2\ell \le q \le k,
    \end{align}
    and
    \begin{align}\label{rangeoindices}
    \ell! \le \sqrt{\norm{u_0}{}} \le \left(k!\right)^{\tfrac{1}{k+1}}, \qquad 2\ell \le q \le k.
\end{align}

Then the complex solution of NSE satisfies
\begin{align*}
    \sup_{t\in(0,\tau_q)} \;\sup_{y\in \mathcal{D}_t}
    \norm{D^qu(\cdot, y, t)}{} + 
    \sup_{t\in(0,\tau_q)} \;\sup_{y\in \mathcal{D}_t}
    \norm{D^qv(\cdot, y, t)}{}
    \,\lesssim & \,
    \norm{D^{q}u_{0}}{} +\\
    & \,\delta_0
    \tau_{k}^{\left(\tfrac{1}{2} + \tfrac{\alpha \beta (\ell + 1)^{1 - \alpha}}{k+1}\right)} \dfrac{(q!)^{\tfrac{1}{q+1}}}{(k!)^{\tfrac{1}{k+1}}} 
    \norm{D^{k}u_{0}}{}^{\tfrac{q+1}{k+1}},
\end{align*}
where \(\tau_q\) is such that
\begin{align*}
    \tau_{q}^{\left(\tfrac{1}{2} +  \tfrac{\alpha \beta (\ell + 1)^{1 - \alpha}}{k + 1}\right)} 
    \leq 
    \left( 1 + \norm{u_0}{2} \left( 1 + \tfrac{1}{2} \log\norm{u_0}{} \right) \right)^{-1} \left( \norm{D^{k}u_{0}}{} \left(q!\right)^{k-q} \right)^{-\tfrac{1}{k+1}},
\end{align*}
and 
\begin{align*}
    \mathcal{D}_{t} = : \{x + iy \in \mathbb{C}^3 : |y| \le ct^{\frac{1}{2}}\}.
\end{align*}
\end{theorem}

\begin{proof}
    To prove Theorem \ref{MainThm1}, it is enough to check that  condition \eqref{cond1} and \eqref{cond2} of Theorem \ref{theorem1} are satisfied. 
    
Since \(\tau_k \in(0,1), \, \beta>0,\,\) \(0<\alpha < 1,\) and  \(q \leq k\), one has
\begin{align*}
    (q+1)^{\alpha} \leq (k+1)^{\alpha} \implies  \dfrac{1}{(q+ 1)^{\alpha}} - \dfrac{1}{(k + 1)^{\alpha}} \geq 0.
\end{align*}
So, for \(q \leq k:\)
    \begin{align*}
        \tau_{k}^{ \beta\left(\tfrac{1}{(q + 1)^{\alpha}} - \tfrac{1}{(k + 1)^{\alpha}}\right)} \leq 1
    \end{align*}
    is satisfied.

    Feeding the bridge inquality into the systems conditions we obtain
\begin{align*}
\sum_{i=\ell}^{q-\ell}
\binom{q}{i}
\mathcal{M}_{i,k}^{\,i+1}
\mathcal{M}_{q-i,k}^{\,q-i+1}
&=
\sum_{i=\ell}^{q-\ell}
\frac{q!}{(q-i)!\,i!}
\left[
\tau_k^{
\beta\left(
\frac{1}{(i+1)^\alpha}
-
\frac{1}{(k+1)^\alpha}
\right)}
\frac{(i!)^{\frac{1}{i+1}}}
     {(k!)^{\frac{1}{k+1}}}
\right]^{i+1}
\\
&\times
\left[
\tau_k^{
\beta\left(
\frac{1}{(q-i+1)^\alpha}
-
\frac{1}{(k+1)^\alpha}
\right)}
\frac{((q-i)!)^{\frac{1}{q-i+1}}}
     {(k!)^{\frac{1}{k+1}}}
\right]^{q-i+1}\\
&\leq  
\dfrac{(q!)^{\tfrac{1}{q+1}}}{(k!)^{\tfrac{1}{k+1}}} (q!)^{\tfrac{k-q}{k+1}} \sum_{i = \ell}^{q - \ell} \tau_k^{\beta\Lambda_{i,q,k}},
\end{align*}
where
\[
\Lambda_{i,q,k}
:=
(i+1)^{1-\alpha}
+(q-i+1)^{1-\alpha}
-\frac{q+2}{(k+1)^\alpha}, \qquad
\ell\le i\le q-\ell.
\]

For the inequality to be valid we need \(\Lambda_{i,q,k} \geq 0.\) Since \(i+1\le k+1\) and \(q-i+1\le k+1\), we have
\[
(i+1)^{-\alpha},\,(q-i+1)^{-\alpha}\ge (k+1)^{-\alpha}.
\]
Hence
\[
(i+1)^{1-\alpha}+(q-i+1)^{1-\alpha}
\ge
\frac{(i+1)+(q-i+1)}{(k+1)^\alpha}
=
\frac{q+2}{(k+1)^\alpha},
\]
and therefore
\[
\Lambda_{i,q,k}\ge0.
\]

Estimating \(\Lambda\) from below: Given that  $2\ell\le q\le k$ and $\ell\le i\le q-\ell$, we have
\( \ell+1\le i+1,\;q-i+1\le k.\)
Hence, for $0<\alpha\le1$,
\begin{align*}
\Lambda_{i,q,k}
&=(i+1)^{1-\alpha}
 +(q-i+1)^{1-\alpha}
 -\frac{q+2}{(k+1)^\alpha} \\
&\ge
2(\ell+1)^{1-\alpha}
\left[
1-\left(\frac{k}{k+1}\right)^\alpha
\right] \\
&\ge
\frac{2\alpha(\ell+1)^{1-\alpha}}{k+1},
\end{align*}
where we used
\[
1-(1-x)^\alpha\ge \alpha x,
\qquad 0\le x\le1.
\]
Since $0<\tau_k<1$ and $\beta>0$, it follows that
\[
\tau_k^{\beta\Lambda_{i,q,k}}
\le
\tau_k^{
\frac{2\alpha\beta(\ell+1)^{1-\alpha}}{k+1}
}.
\]
Consequently,
\[
\sum_{i=\ell}^{q-\ell}
\binom{q}{i}
\mathcal M_{i,k}^{\,i+1}
\mathcal M_{q-i,k}^{\,q-i+1}
\lesssim
\Phi(q,k)\,
\tau_k^{
\frac{2\alpha\beta(\ell+1)^{1-\alpha}}{k+1}
}
\lesssim 1,
\]
provided $\tau_k$ is chosen sufficiently small; in particular, it is
sufficient to take
\[
0<\tau_k<
\min\left\{
1,\,
\Phi(q,k)^{
-\frac{k+1}
{2\alpha\beta(\ell+1)^{1-\alpha}}
}
\right\},
\]
where 
\[\Phi(q, k) =\dfrac{(q!)^{\tfrac{1}{q+1}}}{(k!)^{\tfrac{1}{k+1}}} (q!)^{\tfrac{k-q}{k+1}} \]
For the second condition  \eqref{cond2}, set
\[
a_i:=\frac{(\ell+1)(i+\frac32)}{\ell+\frac32},
\qquad 0\le i\le \ell.
\]
Inserting the bridge into \eqref{cond2}, we obtain
\begin{align*}
    \|u_0\|_2
\sum_{i=0}^{\ell}\binom{q}{i}
\mathcal M_{\ell,k}^{a_i}
\mathcal M_{q-i,k}^{q-i+1}
\left((k!)^{\frac1{k+1}}\|u_0\|\right)^{
\frac{(\frac32-1)(\ell-i)}{\ell+\frac32}}
\lesssim &\,
\frac{\|u_0\|_2\,q! \sqrt{\norm{u_0}{}}}
     {(k!)^{\frac{q+2}{k+1}}}
\sum_{i=0}^{\ell}
\dfrac{1}{i!}
\left(\frac{(\ell!)}{ \sqrt{\norm{u_0}{}}}\right)^{{\tfrac{i+\frac32}{\ell+ \tfrac32}}} \times\\
& \tau_k^{\beta\Gamma_{i,q,k}}, 
\end{align*}
where
\[
\Gamma_{i,q,k}
:=
a_i\left(
\frac1{(\ell+1)^\alpha}
-\frac1{(k+1)^\alpha}
\right)
+
(j-i+1)\left(
\frac1{(q-i+1)^\alpha}
-\frac1{(k+1)^\alpha}
\right) \geq 0.
\]
Dropping the second term and noting that \(a_i \geq a_0 = \tfrac{\tfrac{3}{2}(\ell + 1)}{\ell + \tfrac{3}{2}} \), we have
\begin{align*}
    \Gamma_{i,q,k} 
    &\geq  \left((\ell + 1)^{-\alpha} - (k + 1)^{-\alpha}\right) \dfrac{3(\ell +1)}{2\ell +3} = (\ell + 1)^{1-\alpha}\left[ 1 - \left(\dfrac{\ell + 1}{k + 1}\right)^{\alpha} \right] \dfrac{3}{2\ell + 3}\\
    &\geq (\ell + 1)^{1-\alpha} \dfrac{\alpha (k - \ell)}{k+1} \dfrac{3}{2\ell + 3}\\
    &\geq (\ell + 1)^{1-\alpha} \dfrac{3\alpha \ell}{(k+1)(2\ell + 3)}\\
    &\geq \dfrac{\alpha (\ell + 1)^{1-\alpha}}{k + 1}.
\end{align*}
Here we used the fact that \(1 - r^\alpha \geq \alpha(1 - r)\) (since \(f(r) = r^{\alpha}\) is concave down, for \(0\leq r \leq 1\)), for \(0 < \alpha \leq 1\), and \(k \geq 2\ell \implies k -\ell \geq \ell\) for sufficiently large \(\ell,\) say, \(\ell \geq 3.\) Since $0<\tau_k<1$, we have
\[ \tau_k^{\beta\Gamma_{i,q,k}}
\le
\tau_k^{
\tfrac{\alpha\beta(\ell+1)^{1-\alpha}}{k+1}
}. \]
Summing over $0\le i\le\ell$ and using \eqref{rangeoindices} we obtain the logarithmic bound
\[
\ell\lesssim\frac12\log\|u_0\|.
\]
To this end
\[
\|u_0\|_2
\sum_{i=0}^{\ell}\binom{q}{i}
\mathcal M_{\ell,k}^{a_i}
\mathcal M_{q-i,k}^{q-i+1}
\left((k!)^{\frac1{k+1}}\|u_0\|\right)^{
\frac{(\frac32-1)(\ell-i)}{\ell+\frac32}}
\lesssim
\Psi(q,k)\,
\tau_k^{
\frac{\alpha\beta(\ell+1)^{1-\alpha}}{k+1}
},
\]
where
\[
\Psi(q,k)
:=
\|u_0\|_2
\left(1+\frac12\log\|u_0\|_{\infty}\right)
\frac{(q!)^{\frac1{q+1}}}{(k!)^{\frac1{k+1}}}
(q!)^{\frac{k-q}{k+1}}.
\]

Thus,
\[
\|u_0\|_2
\sum_{i=0}^{\ell}\binom{q}{i}
\mathcal M_{\ell,k}^{a_i}
\mathcal M_{q-i,k}^{q-i+1}
\left((k!)^{\frac1{k+1}}\|u_0\|\right)^{b_i}
\lesssim
{\Psi}(q,k)\,
\tau_k^{\frac{\alpha\beta(\ell+1)^{1-\alpha}}{k+1}}
\lesssim 1,
\]
provided
\[
0<\tau_k<
\min\left\{
1,\,
{\Psi}(q,k)^{
-\frac{k+1}{\alpha\beta(\ell+1)^{1-\alpha}}
}
\right\}.
\]
Theorem \ref{theorem1} applied yields
\begin{align*}
    \norm{D^qu(t)}{} \leq \norm{D^qu_0}{} + \left(\Phi(q,k) + \Psi(q,k)\right)\norm{D^ku_0}{}^{\tfrac{q+2}{k+1}} \left(\tau_{k}\right)^{\left(\tfrac{k+1+2\alpha \beta \left(\ell+1\right)^{1-\alpha}}{2(k+1)}\right)}
\end{align*}
To complete the proof, choose \(\tau_q\) such that
\begin{align*}
    \tau_{q} \leq \left\{\left(\Phi(q,k) + \Psi(q,k)\right)^{-1}
    \norm{D^{k}u_0}{\infty}^{-\tfrac{1}{k+1}}
    \right\}^{\tfrac{2(k+1)}{(k+1)+2\alpha \beta(\ell+1)^{1-\alpha}}}
\end{align*}
\end{proof}

From the proof of Theorem \ref{MainThm1}, we obtain the lower bound on the analyticity radius, which is the square root of the time of existence 
\begin{align}
\label{boss}
    \rho_{k} \geq C\left(\alpha, \beta, \norm{u_0}{2}, \log(\norm{u_0}{\infty})\right)^{-1} \norm{D^{k}u_0}{\infty}^{-\tfrac{1}{(k+1)+2\alpha\beta(\ell+1)^{1-\alpha}}}
\end{align}
Meanwhile, Theorem~\ref{scalebound} gives the characteristic sparseness scale as
\[
r_k
\sim
\|D^k u_0\|_{\infty}^{-\tfrac{1}{{k+\tfrac{3}{2}}}}.
\]

Thus, in order for the analyticity scale to dominate the  scale of sparseness, it is sufficient to require
\[
\dfrac{1}{
k+1+2\alpha\beta(\ell+1)^{1-\alpha}
}
<
\dfrac{1}{k+\tfrac{3}{2}},
\]

which implies that
\[
2\alpha\beta(\ell+1)^{1-\alpha}
>
\dfrac{1}{2}.
\]
Since, \(\alpha\) is arbitrary fixed in \((0, 1)\) and \(\ell + 1 > 1,\) we can choose any \(\beta \geq \tfrac{1}{4\alpha},\) such that the gap is positive. Generally, for any fixed \( \alpha \in (0, 1)\) and \(\beta > 0\), we can choose \(\ell \approx \log\left(\norm{u_0}{\infty}\right)\) sufficiently large such that the required condition
\begin{align*}
    \alpha\beta(\ell+1)^{1-\alpha}>\dfrac{1}{4}
\end{align*}
is satisfied. Indeed, under this condition the exponent gap
\[
\Delta_k
:=
\dfrac{1}{{k+\tfrac{3}{2}}}
-
\frac1{
k+1+2\alpha\beta(\ell+1)^{1-\alpha}
} =  \frac{
2\alpha\beta(\ell+1)^{1-\alpha}-\frac12
}{
\left(k+\frac32\right)
\left(k+1+2\alpha\beta(\ell+1)^{1-\alpha}\right)
} > 0
\]
is strictly positive.

Consequently,
\[
\dfrac{\rho_k}{r_k}
\gtrsim
\|D^k u_0\|_{\infty}^{\Delta_k},
\]
up to the fixed constants entering the analyticity and sparseness
estimates. Hence, in the large-derivative regime,
\[
\rho_k\gtrsim r_k,
\]
which closes the analyticity--sparseness comparison.

\subsection{Proof of Theorem \ref{Mainthm2}}

\begin{proof}[Proof of Theorem \ref{Mainthm2}] Assume, for contradiction, that $T^*$ is the first blow-up time. Set 
\begin{equation}
\label{theta-cond}
    \Theta:=\alpha\beta(\ell+1)^{1-\alpha}  > \frac{1}{4}.
\end{equation}
By hypothesis, there exists \(\varepsilon >0\) and \(K_0 \geq 2\ell\) such that 
\begin{align}
    \label{escape-property}
    \|D^k u(t')\|_{\infty} > \|D^k u(t)\|_{\infty},
    \qquad T^{*} - \varepsilon < t<t'<T^*,
\end{align}
for every integer \(k \geq K_{0}.\) Choose \(t \in (T^{*}-\varepsilon, T^{*})\) sufficiently close to \(T^{*}\), such that \(\ell! \leq \sqrt{\norm{u(t)}{\infty}}\). At this time, select \(k \geq K_0\) large enough such that
\[\norm{u(t)}{\infty} \leq (k!)^{\frac{2}{k+1}} \leq k^2.\]
By Theorem~\ref{MainThm1}, the two nonlinear estimates are controlled
by the common time gain \(\tau_k^{\tfrac{\Theta}{k+1}}.\)
Combining this gain with the usual Navier--Stokes time exponent $1/2$
gives the analyticity exponent
\[
    \frac{1}{2}+\frac{\Theta}{k+1}.
\]
Consider the evolution of the system for
\[
    s\in(t,t+\tau_k)\subset(t,T^*).
\]
For large \(j\) comparable to \(k\) and for large \(\tau_k\), Theorem \ref{MainThm1} gives
\begin{align*}
    \tau_{k} &= \left(1 + \norm{u_0}{2} (1 + \tfrac{1}{2} \log\norm{u_0}{\infty}) 
    \cdot
    \norm{D^ku_0}{\infty}
    \right)^{-\tfrac{2(k+1)}{(k+1) + 2 \alpha \beta (l+1)^{1-\alpha}} \tfrac{1}{k+1}} \\
    & \geq 
    \tfrac{1}{\tilde{c}(\norm{u_0}{2}, \beta, \alpha, \delta_0)} 
    \left[\log\norm{u(t)}{\infty}\norm{D^ku(t)}{\infty}\right]^{-\tfrac{2(k+1)}{(k+1) + 2 \alpha \beta (l+1)^{1-\alpha}} \tfrac{1}{k+1}}.
\end{align*}
Thus the upper bound in Theorem \ref{MainThm1}, at \(j=k\) becomes
\begin{align}
\label{analytic-estimate}
    \sup_{s\in(t,t+\tau_k)} \;\sup_{y\in \mathcal{D}_s}
    \norm{D^ku(\cdot, y, t)}{} + 
    \sup_{s\in(t, t+\tau_k)} \;\sup_{y\in \mathcal{D}_{s}}
    \norm{D^{k}v(\cdot, y, t)}{}
    \,\lesssim & \,
    \left(1 + \delta_0 \left(\tau_k\right)^{\tfrac{1}{2}+\frac{\Theta}{k+1}}\right) \norm{D^{k}u(t)}{},
\end{align} together with the analysis in \eqref{boss}, the solution admits a spatial analytic extension with radius
$\rho_k(s)$ satisfying
\begin{align*}
    \rho_{k}(s) \geq C\left(\alpha, \beta, \norm{u(t)}{2}, \log(\norm{u(t)}{\infty})\right)^{-1} \norm{D^{k}u(t)}{\infty}^{-\tfrac{1}{(k+1)+2\Theta}}.
\end{align*}
Because \(t < s < T^{*}\) and \(k\geq K_{0}\), \eqref{analytic-estimate} gives 
\begin{align*}
    \|D^k u(s)\|_{\infty} > \|D^k u(t)\|_{\infty},
\end{align*} 
and therefore
\[
    \rho_k(s)
    \gtrsim
    \|D^k u(s)\|_{\infty}^{-\tfrac{1}{k+1+2\Theta}}.
\]
On the other hand, the geometric sparseness estimate in Theorem \ref{scalebound} gives the
level-$k$ sparseness scale
\begin{equation}
\label{rk-estimate}
    r_k(s)
    \sim
    \|D^k u(s)\|_{\infty}^{-\tfrac{1}{k+\frac32}}.
\end{equation}
Hence
\[\frac{\rho_k(s)}{r_k(s)}
\gtrsim 
\|D^k u(s)\|_{\infty}^{\Delta_k}, \]
where
\begin{align*}
    \Delta_k
    =
    \frac{1}{k+\frac32}
    -
    \frac{1}{k+1+2\Theta}
    \notag
    =
    \frac{2\Theta-\frac12}
    {\left(k+\frac32\right)(k+1+2\Theta)}>0
\end{align*}
by \eqref{theta-cond}.

Therefore, as $t\uparrow T^*$ and $\|D^k u(s)\|_{\infty}$ becomes sufficiently
large, the positive power $\|D^k u(s)\|_{\infty}^{\Delta_k}$ dominates  the fixed
and logarithmic factors appearing in Theorem \ref{MainThm1}. That is
\begin{align}
\label{rhodominates-r}
    \rho_k(s) \geq \tfrac{1}{c_1 c_2} (\log\norm{u(s)}{\infty})^{-\tfrac{1}{(k+1)+ 2 \alpha \beta(\ell +1)^{1 - \alpha}}} \norm{D^ku(s)}{\infty}^{\Delta_k} r_k(s) 
\end{align}
Consequently, for sufficiently large \(k\), we have \(r_k(s)\le \rho_k(s)\) provided 
\[\norm{D^{k}u(s)}{\infty}^{\Delta_k} \geq c_1 c_2 (\log\norm{u(s)}{\infty})^{\tfrac{1}{(k+1)+ 2 \alpha \beta(\ell +1)^{1 - \alpha}}}.\]
We now use the geometric information at the scale $r_k(s)$.
Choose $0<\varepsilon<1$, and select a derivative component
$D^\zeta u_m$, $|\zeta|=k$, a sign
$\sigma\in\{-1,1\}$, and a point $x_0$ such that
\begin{equation}
\label{approx-max2}
    \sigma D^\zeta u_m(x_0,s)
    \ge (1-\varepsilon) \norm{D^{\zeta}u(s)}{\infty}.
\end{equation}
Fix $0 < \lambda < 1-\varepsilon$ and consider the signed super-level set
\[
    S_\lambda(s)
    :=
    \left\{
        x\in\mathbb R^3:
        \sigma D^\zeta u_m(x,s)>\lambda \norm{D^{\zeta}u(s)}{\infty}
    \right\}.
\]
Then $x_0\in S_\lambda(s)$.

By the one-dimensional sparseness property, there exist a unit vector
$\nu$ and a radius
\[
    0<r\le r_k(s),
\]
such that
\[
    \frac{
        |S_\lambda(s)\cap(x_0-r\nu,x_0+r\nu)|
    }{2r}
    \le\delta
\]
for some fixed $0<\delta<1$.

Consider the good portion of the line segment defined by
\[K :=
    \left\{
        x\in[-r,r]:
        x_0+\nu x\notin S_\lambda(s)
    \right\}.
\]
Then
\begin{equation}
\label{Ksize}
    |K|\ge 2r(1-\delta),
\end{equation}
and, for every $x\in K$,
\begin{equation}
\label{good-set-bound}
    \sigma D^\zeta u_m(x_0+x\nu,s)
    \le \lambda \norm{D^ku(s)}{\infty}.
\end{equation}

Now \eqref{rhodominates-r} implies
\[
    r\le r_k(s)\le\rho_k(s).
\]
Hence, the complex disk $|z|<r$ lies inside the spatial analyticity
region along the complexified direction $\nu$. Therefore,
\[
    f(z)
    :=
    \sigma D^\zeta u_m(x_0+z\nu,s)
\]
is analytic for $|z|<r$.

{ The analytic estimate from \eqref{analytic-estimate} gives
\begin{equation}
\label{analytic-M}
    |f(z)|
    \le
    A_k(t)\norm{D^{k}u(t)}{\infty},
    \qquad |z|<r,
\end{equation}}
where
\[
    A_k(t) := 1 + \delta\left(\tau_k\right)^{\tfrac{1}{2}+\frac{\Theta}{k+1}} \longrightarrow1
    \qquad\text{as }t\uparrow T^*.
\]
Set \( M:=A_k(t)\norm{D^{k}u(t)}{\infty}.\)
By the analytic estimate \eqref{analytic-estimate}, the complex extension of \(D^ku(s)\) is bounded by \(M\) throughout the disk. Restricting this bound to the real domain we get 
\[
    \norm{D^{k}u(s)}{\infty}\le M,
\]
and \eqref{good-set-bound} yields
\[
    f(x)\le\lambda M,
    \qquad x\in K.
\]

Consider the slit domain
\[
    \Omega=D_r\setminus K.
\]
Because $x_0\in S_\lambda(s)$, we have $0\notin K$. That is, assuming otherwise, for \(0 \in K\), we have that \(x_0 + 0\cdot \nu \not \in S_{\lambda(s)},\) implying \(x_0 \not\in S_{\lambda(s)},\) contradiction. Hence,
$0\in\Omega$. Now, let
\[
    h:=\omega(0,\Omega,K)
\]
denote the harmonic measure of $K$ at the origin. By
\eqref{Ksize} and the standard harmonic-measure estimate for subsets
of a diameter,
\[
    h\ge h(\delta)>0,
\]
where $h(\delta)$ depends only on the sparseness parameter.

Since $\Re f$ is harmonic in $\Omega$, the harmonic-measure maximum
principle gives
\begin{align}
    \Re f(0)
    &\le
    h\,\lambda M+(1-h)M
    \notag\\
    &=
    \bigl(1-h(1-\lambda)\bigr)M.
\label{HM-contraction}
\end{align}
On the other hand, by \eqref{approx-max2},
\[
    \Re f(0)
    \ge(1-\varepsilon)\norm{D^ku(s)}{\infty}.
\]
Combining this with \eqref{HM-contraction}, we obtain
\[\norm{D^ku(s)}{\infty}
    \le
    \frac{1-h(1-\lambda)}{1-\varepsilon}
    A_k(t)\norm{D^ku(k)}{\infty}.
\]

Choose $\varepsilon>0$ sufficiently small so that
\[
    \frac{1-h(1-\lambda)}{1-\varepsilon}<1.
\]
Since $A_k(t)\to1$ as $t\uparrow T^*$, we may then choose the escape
time $t$ sufficiently close to $T^*$ so that

\begin{align}\label{contraction-cond}
     \eta
    :=
    \frac{1-h(1-\lambda)}{1-\varepsilon}A_k(t)
    <1.
\end{align}

Therefore
\[
    \norm{D^ku(s)}{\infty}
    \le\eta \norm{D^ku(t)}{\infty}
    <\norm{D^ku(t)}{\infty}.
\]
This contradicts the escape-time property
\eqref{escape-property}, which requires
\[
   \norm{D^ku(s)}{\infty}>\norm{D^ku(t)}{\infty}.
\]
Thus $T^*$ cannot be a blow-up time. Consequently, the solution
continues beyond $T^*$, and the local analyticity result yields a
spatially analytic extension of $u(\cdot,t)$ beyond $T^*$.
\end{proof}


\newpage

\begin{thebibliography}{9}
\bibitem{ConstantinFoias1988}
Constantin, P. and Foias, C. \textit{Navier--Stokes Equations}.
University of Chicago Press, Chicago, 1988.

\bibitem{Kato1984}
Kato, T. Strong $L^p$-solutions of the Navier--Stokes equation
in $\mathbb{R}^m$, with applications to weak solutions.
Mathematische Zeitschrift, 1984, 187(4): 471--480.
doi: 10.1007/BF01174182.

\bibitem{Fefferman2000} Fefferman, C.~L. Existence and smoothness of the Navier--Stokes equation. In \textit{The Millennium Prize Problems}, Clay Mathematics Institute, Cambridge, MA (2000), 57--67.

\bibitem{RobinsonRodrigoSadowski2016}
Robinson, J. C., Rodrigo, J. L., and Sadowski, W. The Three-Dimensional Navier--Stokes Equations: Classical Theory. Cambridge Studies in Advanced Mathematics, Vol. 157. Cambridge University Press, 2016. doi: 10.1017/CBO9781139095143

\bibitem{Tsai2018}
Tsai, T.-P. Lectures on Navier--Stokes Equations. Graduate Studies in Mathematics, Vol. 192. American Mathematical Society, 2018. doi: 10.1090/gsm/192.

\bibitem{FoiasManleyRosaTemam2001}
Foias, C., Manley, O., Rosa, R. and Temam, R.
\textit{Navier--Stokes Equations and Turbulence}.
Cambridge University Press, Cambridge, 2001.

\bibitem{GrujicXu2024}
Gruji\'c, Z. and Xu, L. Asymptotic criticality of the Navier--Stokes regularity problem. \textit{Journal of Mathematical Fluid Mechanics} \textbf{26} (2024), no. 3, 53.

\bibitem{grujic}
Gruji\'c, Z. and  Xu, L.  \textit{Time-global regularity of the Navier--Stokes system with hyper-dissipation: turbulent scenario},
Ann. PDE \textbf{11} (2025), no.~1, Paper No.~9.

\bibitem{Bradshaw}
Bradshaw, Z., Farhat, A. and Gruji\'c, Z. An algebraic reduction of the `scaling gap' in the Navier–Stokes regularity problem. Arch. Ration. Mech. Anal. \textbf{231} (2019), no. 3, 1983–2005.

\bibitem{Phiri}
Phiri, M.P. On Bridging Analyticity and Sparseness in Hyperdissipative Navier–Stokes Systems. J. Math. Fluid Mech. 28, 75 (2026). https://doi.org/10.1007/s00021-026-01054-1

\bibitem{grujicz}
Gruji\'c, Z. The geometric structure of the super-level sets and regularity for 3D Navier-Stokes equations. Indiana Univ. Math. J. \textbf{50} (2001), no.3, 1309–1317.

\bibitem{Solynin1999}
Solynin, A.Y. Ordering of sets, hyperbolic metrics, and harmonic measures. J Math Sci 95, 2256–2266 (1999). https://doi.org/10.1007/BF02172470

\bibitem{GrujicKukavica1998}
Gruji\'c, Z. and Kukavica, I. Space analyticity for the Navier--Stokes and related equations with initial data in $L^p$. \textit{Journal of Functional Analysis} \textbf{152} (1998), no. 2, 447--466.

\bibitem{farhat}
Farhat, A. and Gruji\'c, Z. On persistence of spatial analyticity in the hyper-dissipative Navier--Stokes models. Indiana Univ. Math. J. (to appear) https://doi.org/10.48550/arXiv.2312.14320

\bibitem{grujic2}
Gruji\'c, Z. A geometric measure-type regularity criterion for solutions to the 3D Navier-Stokes equations. Nonlinearity \textbf{26} (2012), no.1, 289.

\bibitem{Nirenberg}
Nirenberg, L. On elliptic partial differential equations. Ann. Scuola Norm. Sup. Pisa (3) \textbf{13} (1959), 115–162.

\bibitem{Gagliardo}
Gagliardo, E. Ulteriori propriet\`a di alcune classi di funzioni in pi\`u variabili. Ricerche Mat. \textbf{8} (1959), 24–51.

\bibitem{GrujicMohebujjaman2026}
Gruji\'c, Z. and Mohebujjaman, M. On higher-order derivative ratios in turbulent flows. \textit{arXiv preprint arXiv:2605.21501} (2026).

\bibitem{Ransford}
Ransford, T. \textit{Potential Theory in the Complex Plane}. Cambridge University Press, Cambridge, 1995.




\bibitem{AlbrittonBradshaw}
Albritton, D.  and Bradshaw, Z. Remarks on sparseness and regularity of Navier--Stokes solutions, \textit{Nonlinearity} \textbf{35} (2022), no.~6, 2858--2877.

\bibitem{Guberovic2010} Guberovi\'c, R. Smoothness of Koch-Tataru solutions to the Navier-Stokes equations revisited. Discrete and Continuous Dynamical Systems, 2010, 27(1): 231-236. doi: 10.3934/dcds.2010.27.231

\bibitem{Nirenberg1959}
Nirenberg, L. On elliptic partial differential equations.
Annali della Scuola Normale Superiore di Pisa,
Classe di Scienze, Serie III, 1959, 13(2): 115--162.
\end{thebibliography}
\end{document}